\documentclass[12pt]{article}
\usepackage[dvips]{graphicx}
\usepackage{amsmath, amssymb}
\usepackage{amsthm}

\newtheorem{theorem}{Theorem}[section]
\newtheorem{lemma}{Lemma}[section]

\newtheorem{remark}{Remark}[section]
\newtheorem{proposition}{Proposition}[section]

\newcommand{\cK}{{\cal K}}
\newcommand{\cH}{{\cal H}}
\newcommand{\cP}{{\cal P}}
\newcommand{\cQ}{{\cal Q}}
\newcommand{\bR}{{\mathbb R}}

\newlength{\IndentI}
\newlength{\IndentII}
\newlength{\IndentIII}
\newlength{\WidthI}
\newlength{\WidthII}
\newlength{\WidthIII}
\title{Game-theoretic versions of strong law of large numbers for
  unbounded variables}

\author{
  Masayuki Kumon\\
  Risk Analysis Research Center\\
  Institute of Statistical Mathematics\\
  Akimichi Takemura\\
  Graduate School of Information Science and Technology\\
  University of Tokyo\\
  and\\
  Kei Takeuchi\\
  Faculty of International Studies\\
  Meiji Gakuin University }

\date{March, 2006}

\begin{document}
\maketitle

\begin{abstract}
  We consider strong law of large numbers (SLLN) in the framework of
  game-theoretic probability of Shafer and Vovk (2001).  We prove
  several versions of SLLN for the case that Reality's moves are
  unbounded.  Our game-theoretic versions of SLLN largely
  correspond to standard measure-theoretic results.  However
  game-theoretic proofs are different from measure-theoretic
  ones in the explicit consideration of various hedges.  In
  measure-theoretic proofs existence of moments are assumed, whereas
  in our game-theoretic proofs we assume availability of various hedges
  to Skeptic for finite prices.
 \end{abstract}

\noindent
{\it Keywords and phrases:} \ 
Borel-Cantelli lemma,
call option,
Doob's upcrossing lemma,
Kronecker's lemma,
Marcinkiewicz-Zygmund strong law,
martingale convergence theorem.
\section{Introduction}
\label{sec:intro}

In the framework of game-theoretic probability, proof of SLLN is
simple if Reality's moves are bounded. In \cite{kumon-takemura} we
showed that a single simple strategy based on past averages of
Reality's moves forces SLLN for the case of bounded Reality's moves.
For the special case of the coin-tossing game path behavior and
convergence rate of SLLN can be very explicitly stated
(\cite{ktt1},\cite{takemura-suzuki}).  However when Reality's moves
are not bounded, the proof becomes more complicated due to
consideration of availability of hedges to Skeptic.  Under the
requirement of the collateral duty that Skeptic has to keep his
capital always nonnegative, he has to use some form of hedge at each
round.  In Chapter 4 of Shafer and Vovk (2001), Kolmogorov's SLLN is
proved under the availability of the variance hedge (quadratic hedge).
Shafer and Vovk consider the case that the price of the variance hedge
is announced by Forecaster for each round, but for simplicity in this
paper we omit Forecaster from the protocol and consider the case that
hedges carry constant prices throughout the game.  Availability of the
quadratic hedge is natural and convenient.  However the purpose of
this paper is to investigate SLLN under other types of hedges.

In measure-theoretic probability, the usual and most elegant form of
SLLN is stated for the sample average $\bar x_n=(1/n)(x_1 + \dots +
x_n)$ of i.i.d.\ random variables, where only the existence of the
measure-theoretic expected value $E|x_n|<\infty$ is assumed.  However Kolmogorov's SLLN
proved in Chapter 4 of Shafer and Vovk (2001) does not correspond to
this version and a question remains whether a corresponding
game-theoretic result holds or not.  
Some considerations of this problem are given in Chapter 4 of \cite{takeuchi:2004}.
The usual measure-theoretic result
depends strongly on the assumption of identical distribution of the
random variables.  On the other hand the basic feature of the
game-theoretic probability is that the game is a martingale and there
is a question of how to impose identical behavior to Reality at each
round.  In this paper we argue that the assumption of the identical
distribution in measure-theoretic framework can be replaced by the
availability of countable number of weak hedges.

For the most part we follow the standard proofs of SLLN in
measure-theoretic probability.  For example we use truncation and
Kronecker's lemma.  However our proofs differ from standard
measure-theoretic proofs in explicit construction of Skeptic's
strategy which requires Skeptic to observe his collateral duty.  In
addition our proof is more an extension of the proof for the bounded
case of Chapter 3 of Shafer and Vovk (2001), rather than an extension
of their proof in Chapter 4 using the quadratic hedge.

The organization of this paper is as follows.
In Section \ref{sec:notation} we set up notations and give some
preliminary results.  In Section \ref{sec:single-hedge} we prove a
version of SLLN under the assumption of availability of a single hedge.
In Section \ref{sec:with-truncation} we prove a game-theoretic version
of SLLN for i.i.d.\ variables under the assumption of availability of
countable hedges.  We extend it to a Marcinkiewicz-Zygmund strong law
in Section \ref{sec:MZ}.  Finally in Section \ref{sec:discussions} we discuss
various aspects of our proofs and the assumption of availability of
infinite number of hedges.

\section{Notation and preliminaries}
\label{sec:notation}

In this section we summarize our notations and some preliminary
results.  We follow the notation of Shafer and Vovk (2001).  $\xi=x_1
x_2 \dots$ denotes an infinite path of Reality's moves and $\xi^n = x_1 \dots
x_n$ denotes the partial path up to round $n$.  For a strategy ${\cal
  P}$ of Skeptic, $\cK^\cP_n(\xi)=\cK^\cP_n(\xi^n)$ denotes the
capital process.  Starting with a positive initial capital of ${\cal
  K}_0 = \delta  > 0$,
Skeptic observes his collateral duty by using $\cP$ if
\begin{equation}
\label{eq:duty}
\cK^\cP_n(\xi)  \ge 0, \quad \forall\xi, \forall n .
\end{equation}
We also say that ${\cal P}$ satisfies the collateral duty with the
initial capital  $\delta$.  Note that $\cP$ satisfies the collateral duty with initial
capital $\delta$ if and only if 
$\cP/\delta$ satisfies the duty with the initial capital 1.  In view
of this fact, we simply say that $\cP$ satisfies the collateral duty if 
$\cP$ satisfies the duty with some initial capital $\delta > 0$.
When ${\cal P}$ satisfies the collateral duty, the capital process
$\cK^\cP$ is called a (game-theoretic) non-negative martingale.

We call a function $h(x)$ of Reality's move $x$ a {\it
  hedge} if it is non-negative ($h(x) \ge 0, \forall x\in {\mathbb
  R}$) and has a finite price $0 < \nu < \infty$.  Skeptic is allowed
to buy arbitrary amount of $h(x)$ with the unit price $\nu$.  In
Chapter 4 of Shafer and Vovk (2001), they consider the variance hedge
$h(x)=x^2$.  In view of the unbounded forecasting game in Chapter 4 of
Shafer and Vovk (2001), we first consider the following protocol with a
single hedge.

\medskip
\noindent
\textsc{The Unbounded Forecasting Game with a Single Hedge}\\
\textbf{Protocol:}

\parshape=6
\IndentI   \WidthI
\IndentI   \WidthI
\IndentII  \WidthII
\IndentII  \WidthII
\IndentII  \WidthII
\IndentI   \WidthI
\noindent
$\cK_0 :=1$.\\
FOR  $n=1, 2, \dots$:\\
  Skeptic announces $M_n\in\bR$, $V_n \ge 0$.\\
  Reality announces $x_n\in\bR$.\\
  $\cK_n := \cK_{n-1} + M_n x_n  + V_n (h(x_n) - \nu)$\\
END FOR

\medskip
Availability of the variance hedge $h(x)=x^2$ is very convenient, because Skeptic
can then construct a martingale which is a quadratic form of Reality's
moves.  This fact is used by Shafer and Vovk in their proof.
However SLLN can be proved under other hedges.
In Section \ref{sec:single-hedge} we will prove 
that SLLN is forced if the absolute moment hedge of order $1+\epsilon$, $\epsilon>0$,
\[
h(x)= |x|^{1+\epsilon}
\]
is available to Skeptic.  Naturally we are tempted to consider the
absolute moment hedge
\[
h(x)=|x|
\]
in the above protocol, corresponding to the measure-theoretic SLLN of
i.i.d.\ random variables with finite expectation.  However it is
essential to point out that SLLN is not forced under the availability
of $h(x)=|x|$ alone.  Since this fact is important, we state it as a
proposition. The following proposition is stated in view of 
the Marcinkiewicz-Zygmund strong law in Section 
\ref{sec:MZ}.

\begin{proposition}
\label{prop:1}
Consider the unbounded forecasting game with a single hedge
$h(x)=|x|^r$, $r > 0$. There exists no strategy $\cP$ of Skeptic
satisfying the collateral duty, such that $\lim_n \cK^\cP_n=\infty$
whenever $(x_1 + \cdots + x_n)/n^{1/r} \not\rightarrow 0$.
\end{proposition}
Proof of this proposition, following Section 4.3 of Shafer and
Vovk (2001), is given in Appendix \ref{sec:app1}.  Unfortunately it
requires a measure-theoretic argument.

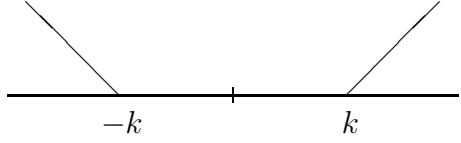
\begin{figure}[htbp]
\setlength{\unitlength}{5mm}
\begin{center}
\begin{picture}(20,4)(-10,0)
\put(-6,0){\line(1,0){12}}
\put(0,-0.2){\line(0,1){0.4}}
\put(-3,0){\line(-1,1){2.5}}
\put(-3.5,-1){$-k$}
\put(3,0){\line(1,1){2.5}}
\put(2.9,-1){$k$}
\end{picture}
\end{center}
\caption{Symmetric call option type hedge}
\label{fig:1}
\end{figure}
Because of Proposition \ref{prop:1} with $r = 1$, we need to assume
that more hedges in addition to $h(x)=|x|$ are available to Skeptic in
order to prove SLLN corresponding to the sample average of
i.i.d.\ random
variables with  finite measure-theoretic expected value $E|x_n|<\infty$.  Let
\[
\cH = \{ h_\lambda \mid \lambda \in \Lambda\}
\]
denote a set of hedges available to Skeptic in each round.
For example in Section  \ref{sec:with-truncation} we consider the set
of symmetric call option type hedges (``strangle hedges'', Chapter 10
of \cite{hull})
\[
\cH = \{ (|x|-k)_+ \mid  k=0,1,2,\dots\},
\]
where $x_+ = \max(0,x)$.  $(|x|-k)_+$ is depicted in Figure
\ref{fig:1}.
We assume that $h_\lambda$ is available to
Skeptic with a constant finite positive price $v_{h_\lambda}$.  Skeptic is allowed
to buy any amount of countable number of hedges $h_1,
h_2,\dots$ from $\cH$.  
If Skeptic buys $V_i \in \bR$ units of $h_i$,
$i=1,2,\dots$, then he is required that the sum 
$\sum_{i=1}^\infty V_i v_{h_i}$ converges to a finite value.
Note that here for a set of hedges we are allowing Skeptic to sell a
hedge ($V_i < 0$), whereas in the case of a single hedge Skeptic can
obviously only buy the hedge.
By allowing Skeptic to sell hedges, he can combine various hedges to
construct a variety of hedges (Chapter 10 of \cite{hull}, Section 9.3 of 
\cite{capinski-zastawniak}).
Based on these considerations we set up the following protocol.

\medskip
\noindent
\textsc{The Unbounded Forecasting Game with a Set of Hedges}\\
\textbf{Protocol:}

\parshape=7
\IndentI   \WidthI
\IndentI   \WidthI
\IndentII  \WidthII
\IndentIII  \WidthIII
\IndentII  \WidthII
\IndentII  \WidthII
\IndentI   \WidthI
\noindent
$\cK_0 :=1$.\\
FOR  $n=1, 2, \dots$:\\
  Skeptic announces $M_n\in\bR$, $h_{n1}, h_{n2}, \dots \in \cH$,
  $V_{n1}, V_{n2}, \dots \in \bR$\\
 s.t.\ $\sum_i V_{ni} v_{h_{ni}}$ converges to a finite value.\\
  Reality announces $x_n\in\bR$.\\
  $\cK_n := \cK_{n-1} + M_n x_n 
   + \sum_i V_{ni} (h_{ni}(x_n)-v_{h_{ni}})$.\\
END FOR

\medskip

In our proofs we combine Skeptic's strategies to force intersection of
events. {}From Section 3.2 of Shafer and Vovk (2001), a strategy $\cP$
weakly forces an event $E$ if it satisfies the collateral duty and
$\limsup_n \sum_n \cK^\cP_n(\xi)=\infty$ for every $\xi \not\in E$.
In this case we also say that $E$ happens almost surely.  If
$\limsup_n$ is replaced by $\lim_n$, then $\cP$ forces $E$.  Now
consider two events $E_1$ and $E_2$.  We say that a strategy $\cP$
weakly forces $E_2$ {\em conditional on} $E_1$ if it satisfies the
collateral duty and
\[
\limsup_n \cK^{\cQ}_n(\xi) = \infty, \qquad \forall \xi \in E_1 \cap E_2^C.
\]

Now we state the following lemma, which is slightly stronger than
Lemma 3.2 of Shafer and Vovk (2001).

\begin{lemma}
\label{lem:conditional}
Suppose that Skeptic can weakly force $E_1$ and furthermore he can
weakly force $E_2$ conditional on $E_1$.  Then he can weakly 
force $E_1 \cap E_2$.
\end{lemma}

\begin{proof}
Let $\cP_1$ denote a strategy weakly forcing $E_1$ and let $\cP_2$
denote a strategy weakly forcing $E_2$ conditional on $E_1$.  Let 
$\cP = (1/2) (\cP_1 + \cP_2)$. Note that
\[
(E_1 \cap E_2)^C = E_1^C \cup ( E_1\cap E_2^C).
\]
For $\xi \in E_1^C$ $\limsup_n \cK^\cP_n (\xi)=\infty$ since
$\limsup_n \cK^{\cP_1}_n(\xi)=\infty$.  Similarly for 
For $\xi \in E_1\cap E_2^C$ $\limsup_n \cK^\cP_n (\xi)=\infty$ since
$\limsup_n \cK^{\cP_2}_n(\xi)=\infty$.
\end{proof}

It is clear that Lemma \ref{lem:conditional} can be generalized to the 
sequence of events $E_1, E_2, \dots$, such that $E_i$ is weakly forced
conditional on $E_1 \cap \dots \cap E_{i-1}$.

Finally we state and discuss the game-theoretic martingale convergence
theorem given in Lemma 4.5 of Shafer and Vovk (2001).  
\begin{lemma}
\label{lem:upcrossing}
 A non-negative martingale $\cK^\cP$ converges to a
  non-negative finite value almost surely.
\end{lemma}

As seen from the proof of Lemma 4.5 of Shafer and Vovk (2001) this
theorem is based on Doob's upcrossing lemma in the game-theoretic
setting.  We use this lemma in our proofs in an essential way.  As
discussed at the beginning of this section, when we say that $\cK^\cP$
is a non-negative martingale, it means that Skeptic 
observes his collateral duty (\ref{eq:duty}) with the strategy $\cP$
starting with a positive initial capital $\cK_0 > 0$.
In this case he can construct another strategy $\cQ$ satisfying the
collateral duty starting with an arbitrary small initial capital $\delta > 0$
such that
\[
\lim_n \cK^{\cQ}_n(\xi) = \infty
\]
whenever $\cK^\cP_n$ does not converge.  As in Section 4.2 of Shafer
and Vovk (2001) or Chapter 12 of Williams (1991) we use 
Lemma \ref{lem:upcrossing} in conjunction with Kronecker's lemma.

\section{SLLN with a single hedge}
\label{sec:single-hedge}

In this section we give sufficient conditions for SLLN in the
unbounded forecasting game with a single hedge.  For simplicity we
only consider symmetric hedge $h(x)=h(|x|)$ depending only on $|x|$.  We
assume several conditions for $h(|x|)\ge 0$.

\begin{align}
&\text{(A1)\qquad For some $c>0$,}\qquad      h(|x|) \ge |x|  \text{
  for }\ |x|\ge c. \\
&\text{(A2)\qquad For some $c>0$ and for all $\alpha\ge 1$} \nonumber \\
& \qquad \qquad\qquad 
 \frac{h(|x|)}{|x|^\alpha} \ \text{is monotone increasing or
   decreasing for}\ |x|\ge c.\\
&\text{(A3)\qquad For some $c >0$,} \qquad \sum_{n > c}^\infty
\frac{1}{h(n)} < \infty. 
\label{eq:sum1}
\end{align}

In our proof the condition (A3) is essential for SLLN with a single
hedge, as shown in Proposition \ref{prop:a3} below.  On the
other hand (A2) and the symmetry of $h$ are assumed for convenience
for our proofs.  $c>0$ in the conditions can be easily handled and for
simplicity we assume $c=0$ in our proofs below.  By (A2), there exists
some $\alpha_0>0$ such that $h(|x|)$ is monotone increasing in $|x|$ for
$\alpha>\alpha_0$ and monotone decreasing in $|x|$ for
$\alpha<\alpha_0$.

Now we state the following theorem.
\begin{theorem}
\label{thm:1}
Suppose that a single hedge $h(x)$ satisfying {\rm (A1)--(A3)}
is available to Skeptic.  Then 
in the unbounded forecasting game with the single hedge $h(x)$, 
Skeptic can force $\bar x_n \rightarrow 0$.
\end{theorem}

Take $h(x)=|x|^{1+\epsilon}$, $\epsilon > 0$,
then (A1)--(A3)
hold and SLLN is forced.  SLLN is forced even for
\[
h(x)= |x| (\log |x|)^2.
\]
However as shown in Proposition \ref{prop:1}, SLLN is not forced for
$h(x)=|x|$.

Before starting the proof of Theorem \ref{thm:1} we show that the
condition (A3) is also necessary for the existence of a strategy
weakly forcing SLLN.  

\begin{proposition}
\label{prop:a3}
Consider  $h(x)\ge 0$  with $h(0)=0$ and $\sum_n 1/h(n) =
\infty$. Then  in the unbounded forecasting game with this single hedge $h(x)$, 
there exists no strategy $\cP$ of Skeptic
satisfying the collateral duty, such that $\lim_n \cK^\cP_n=\infty$
whenever $(x_1 + \cdots + x_n)/n \not\rightarrow 0$.
\end{proposition}

Proof of this proposition is given in Appendix \ref{sec:app1}.

The rest of this section is devoted to a proof of Theorem \ref{thm:1}
in a series of lemmas.  By Lemma 3.1 of Shafer and Vovk (2001) we only need
to show that Skeptic can weakly force $\bar x_n \rightarrow 0$.

\begin{lemma} 
\label{lem:BC0}
Let 
\[
E_1 = \{ \xi \mid \sum_n \frac{h(x_n)}{h(n)} < \infty \}.
\]
Under the conditions {\rm (A1)--(A3)}
Skeptic can force $E_1$.
\end{lemma}

\begin{proof}
By (A3)  let $C=\sum_{n} 1/h(n) < \infty$. Consider the following strategy $\cP$
\[
M_n\equiv 0, \qquad  V_n = \frac{1}{C \nu h(n)}.
\]
where $0 < \nu < \infty$ is the price of the hedge $h$. For this strategy, starting
with the initial capital of $\cK_0=1$,  the capital
process $\cK_n$ is written as
\begin{align*}
\cK_n &= 1 + \sum_{i=1}^n \frac{1}{C \nu h(i)} (h(x_i)-\nu)\\
&=  1 - \frac{1}{C} \sum_{i=1}^n \frac{1}{h(i)} 
+ \frac{1}{C\nu} \sum_{i=1}^n \frac{h(x_i)}{h(i)} \\
& \ge \frac{1}{C\nu} \sum_{i=1}^n \frac{h(x_i)}{h(i)}.
\end{align*}
Therefore $\cP$ satisfies the collateral duty and on 
$E_1^C$ $\cK_n$ diverges to $+\infty$.  Therefore $\cP$
forces $E_1$.  
\end{proof}

Note that the same argument with $C=\sum_n 1/n^2$ shows that Skeptic
can force
\begin{equation}
\label{eq:E1'}
E_1' = \{ \xi \mid \sum_n \frac{h(x_n)}{n^2} < \infty \}.
\end{equation}
Furthermore Lemma \ref{lem:BC0} implies the following Borel-Cantelli
type result.

\begin{lemma} 
\label{lem:BC1}
Let
\begin{equation}
\label{eq:E2}
E_2 = \{ \xi \mid |x_n| \ge n  \text{\rm \ for only finite number of } n\}.
\end{equation}
Under the conditions {\rm (A1)--(A3)} 
Skeptic can force $E_2$.
\end{lemma}

\begin{proof}
By (A2) $h(|x|)/|x|$ is monotone.  If it is monotone decreasing (A3)
can not hold.  Therefore $h(|x|)/|x|$ has to be monotone increasing 
and $h(|x|)$ is itself monotone increasing. Therefore for  $z>0$ 
\[
\frac{h(z)}{h(n)} \ge I_{[n,\infty)}(z), 
\]
where $I_{[n,\infty)}(\cdot)$ is the indicator function of the
interval $[n,\infty)$.  It follows that $E_1 \subset E_2$. 
\end{proof}

It should be noted that this lemma is essentially the first part of
Borel-Cantelli lemma.  For convenience we state a game-theoretic
version of the first part of Borel-Cantelli lemma.  The proof is the
same as in Lemma \ref{lem:BC0} and omitted.

\begin{lemma}
\label{lem:BC1-general} {\rm (The first part of Borel-Cantelli)} \quad 
  Let  $E_1, E_2,
  \dots$ be a sequence of events such that the sum of the upper 
  probabilities is finite $\sum_n \bar P(E_n) <  \infty$. Then Skeptic can force 
\[
(\limsup_n E_n)^C = \{ E_n \ \text{\rm only for finite } n \}.
\]
\end{lemma}

The following lemma concerns the evaluation of the variance of
truncated variables in the usual proof of SLLN.

\begin{lemma}
\label{lem:trunc}
Let
\begin{equation}
\label{eq:E3}
E_3 = \{ \xi \mid \sum_n \frac{x_n^2}{n^2}I_{\{|x_n|\le n\}} < \infty \}.
\end{equation}
Under the conditions {\rm (A1)--(A3)} 
Skeptic can force $E_3$.
\end{lemma}

\begin{proof}
First consider the case that $h(x)/x^2$ is monotone increasing.
Then adjusting some constants we can assume $h(x)\ge
x^2$ for all $x$ without loss of generality.  Then
\[
\sum_n \frac{x_n^2}{n^2}I_{\{|x_n|\le n\}} \le 
\sum_n \frac{x_n^2}{n^2} \le \sum_n \frac{h(x_n)}{n^2}
\]
and  $E_1' \subset E_3$, where $E_1'$ is given in 
(\ref{eq:E1'}).  Therefore Skeptic can force $E_3$. 

Next consider the case that $h(x)/x^2$ is monotone decreasing. For $0
< z \le n$ we have
\[
\frac{h(z)}{z^2} \ge \frac{h(n)}{n^2}.
\]
Multiplying both sides by $n^2/h(z)$ we have
\[
\frac{z^2}{n^2} \le \frac{h(z)}{h(n)}.
\]
Then
\[
\sum_n \frac{x_n^2}{n^2}I_{\{|x_n|\le n\}} \le 
\sum_n \frac{h(x_n)}{h(n)}
\]
and $E_1 \subset E_3$.  
\end{proof}

From Lemma \ref{lem:BC1} and Lemma \ref{lem:trunc}  Skeptic can force
$E_2 \cap E_3$.

\begin{lemma}  \label{lem:bound1} 
Let $ 0 < c \le 1/[2(1+\nu/h(1))]$. Then for all $x$
\[
 -c \frac{|x|}{n} + c \frac{ h(x)-\nu}{h(n)} \ge -\frac{1}{2}.
\]
\end{lemma}

\begin{proof}  Since $h(z)/z$ is  increasing in $z>0$, for $z \ge n$
  we have $h(n)/n \le h(z)/z$.  Multiplying by $z/h(n)$ we have
\[
\frac{h(z)}{h(n)} - \frac{z}{n} \ge 0, \quad z \ge n.
\]
For $0 \le z \le n$ obviously 
\[
\frac{h(z)}{h(n)} - \frac{z}{n} \ge -1.
\]
Therefore for all $z \ge 0$ we have
\[
\frac{h(z)-\nu}{h(n)} - \frac{z}{n} \ge -1 - \frac{\nu}{h(n)} \ge 
-1 - \frac{\nu}{h(1)}
\]
and this proves the lemma.
\end{proof}

Finally the following lemma proves Theorem
\ref{thm:1} by Kronecker's lemma.
\begin{lemma}  
\label{lem:e4}
Let
\begin{equation}
\label{eq:E4}
E_4 = \{ \xi \mid  \sum_n \frac{x_n}{n} \ \text{\rm converges to a finite
 value} \}.
\end{equation}
Under the conditions {\rm (A1)--(A3)}
Skeptic can weakly force $E_4$ conditional on $E_1$.
\end{lemma}

\begin{proof} Let
$ 0 < \epsilon \le 1/[2(1+\nu/h(1))]$.
Consider the following strategy $\cP^+$:
\[
M_n = \epsilon \cK_{n-1} \frac{1}{n}, \quad V_n = \epsilon \cK_{n-1}
\frac{1}{h(n)}.
\]
Then by Lemma \ref{lem:bound1} 
\[
\cK_n = \cK_{n-1} (1 +  \epsilon \frac{x_n}{n} + \epsilon \frac{h(x_n) - \nu}{h(n)})
 \ge \frac{1}{2} \cK_{n-1}
\]
and $\cP^+$ satisfies the collateral duty.  Similarly the strategy
$\cP^-$ with $M_n = -\epsilon \cK_{n-1}/n$, $V_n = \epsilon
\cK_{n-1}/h(n)$ satisfies the collateral duty.  By the game-theoretic
martingale convergence theorem (Lemma \ref{lem:upcrossing}) 
both $\cK^{\cP^+}_n$ and  $\cK^{\cP^-}_n$ converge to a non-negative
finite value almost surely. Then 
both $\log
\cK^{\cP^+}_n$ and $\log \cK^{\cP^-}_n$ converge to a finite value
or $-\infty$ almost surely.

As in Lemma 3.3 of Shafer and Vovk (2001) we use the inequality
$t \ge \log (1+t) \ge t- t^2$ for all $t\ge -1/2$.  Then the logarithm
of the capital
process for $\cP^+$ starting with $\cK_0=1$ is bounded as
\begin{align}
\epsilon \sum_{i=1}^n \big(\frac{x_i}{i} - 
\frac{h(x_i)-\nu}{h(i)}\big)&\ge
\log \cK^{\cP^+}_n  \nonumber \\
&\ge  \epsilon \sum_{i=1}^n \big(\frac{x_i}{i} -
\frac{h(x_i)-\nu}{h(i)}\big)
 - \epsilon^2\sum_{i=1}^n \big(\frac{x_i}{i} - 
\frac{h(x_i)-\nu}{h(i)}\big)^2.
\label{eq:bound1}
\end{align}
On $E_1$, each of the following infinite sums is finite.
\[
\sum_n \frac{h(x_n)}{h(n)}, \quad
\sum_n \frac{\nu}{h(n)}, \quad
\sum_n \frac{x_n^2}{n^2}, \quad
\sum_n \frac{h(x_n)^2}{h(n)^2}, \quad
\sum_n \frac{\nu^2}{h(n)^2}.
\]
By the inequality 
\[
(a_1+\dots+a_m)^2 \le m (a_1^2 + \dots + a_m^2)
\]
on $E_1$ the second term on the right-hand side of 
(\ref{eq:bound1}) converges to a finite value:
\[
\sum_{n=1}^\infty \big(\frac{x_n}{n} - 
\frac{h(x_n)-\nu}{h(n)}\big)^2  < \infty. 
\]
Therefore conditional on $E_1$ \  $\cP^+$ weakly forces
\[
\limsup_n \sum_{i=1}^n \frac{x_i}{i} < \infty.
\]
Similarly conditional on $E_1$ \  $\cP^-$ weakly forces
\[
\liminf_n \sum_{i=1}^n \frac{x_i}{i} < -\infty.
\]
It follows that the case 
$\lim \log \cK^{\cP^+}_n=-\infty$ is eliminated and
$\log \cK^{\cP^+}_n$
converges to a finite value almost surely.

Now consider (\ref{eq:bound1}) for the interval $n \le i \le n'$.
Then
\begin{align*}
\epsilon \sum_{i=n}^{n'} \big(\frac{x_i}{i} - 
\frac{h(x_i)-\nu}{h(i)}\big)
&\ge
\log \cK^{\cP^+}_{n'} - \log \cK^{\cP^+}_{n-1}  \nonumber \\
&\ge  \epsilon \sum_{i=n}^{n'} \big(\frac{x_i}{i} - 
\frac{h(x_i)-\nu}{h(i)}\big)
 - \epsilon^2\sum_{i=n}^{n'} \big(\frac{x_i}{i} - 
\frac{h(x_i)-\nu}{h(i)}\big)^2.
\end{align*}
Now by Cauchy criterion we see that $\sum_n x_n/n$ converges almost
surely.
\end{proof}

\section{SLLN with countable hedges}
\label{sec:with-truncation}

In this section we prove a version of game-theoretic SLLN which
corresponds to the usual measure-theoretic SLLN for the sample average
of i.i.d.\ random variables with finite expectation.  As shown in
Proposition \ref{prop:1}, the availability of a single $h(x)=|x|$ is
not sufficient.  It seems that an essential ingredient of
measure-theoretic proofs of SLLN for this case is that the expected
values of truncation are uniformly bounded by the assumption of
identical distribution.  Hence we consider that countable number
of hedges are available with constant prices at each round of the game.
We assume that the prices are given in such a way that the game is
coherent, i.e.\ the game does not present an arbitrage opportunity to
Skeptic (see Section 7.1 of \cite{shafer/vovk:2001}
or \cite{takemura-suzuki}).

As mentioned in Section \ref{sec:notation}, for our game-theoretic
version of SLLN we assume that the set of symmetric call option type hedges
with integral exercise prices $k=0,1,2,\dots$
\begin{equation}
\label{eq:countable-hedges}
\cH = \{ h_k(x)= (|x|-k)_+ \mid  k=0,1,2,\dots\}
\end{equation}
are available to Skeptic.  
In particular $|x|=(|x|-0)_+$ is available to Skeptic.  
Let $\nu_k$ denote the price of $h_k(x)=(|x|-k)_+$.
We also assume that Skeptic is allowed to
sell hedges and combine them, as long as he observes
his collateral duty.  For example
he can create a new hedge
\[
(|x| - k)_+ - (|x| - k-1)_+  =
\begin{cases}
0,&   |x| \le k\\
|x| - k, & k < |x|\le  k+1\\
1, & |x| > k+1.
\end{cases}
\]
This new hedge carries the price of $\nu_k - \nu_{k+1} \ge 0$.  We
may call this hedge ``symmetric bull spread'' (c.f.\ Chapter 10 of \cite{hull}).

\begin{figure}[htbp]
\setlength{\unitlength}{5mm}
\begin{center}
\begin{picture}(16,3)(-1,0)
\put(-1,0){\line(1,0){14}}
\put(0,-0.2){\line(0,1){0.4}}
\put(4,0){\line(1,1){2}}
\put(3.1,-1){$k-1$}
\put(6,2){\line(1,0){2}}
\put(8,2){\line(1,-1){2}}
\put(5.8,-1){$k$}
\put(6.9,-1){$k+1$}
\put(9.4,-1){$k+2$}
\put(6,-0.1){\line(0,1){0.2}}
\put(8,-0.1){\line(0,1){0.2}}
\end{picture}
\end{center}
\caption{Symmetric trapezoidal hedge (positive side only)}
\label{fig:3}
\end{figure}
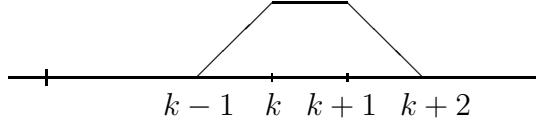

For truncation arguments below we also consider ``symmetric
trapezoidal hedge''.
For $k\ge 1$ define
\begin{align*}
T_k(x) &= (|x| - (k-1))_+ - (|x| - k)_+  - \big(
(|x| - (k+1))_+ - (|x| - (k+2))_+ \big)\\
& =
\begin{cases}
0,&   |x| \le k-1\\
|x| - (k-1), & k-1 < |x| \le  k\\
1, & k < |x| \le k+1\\
k+2 - |x|, & k+1 < |x| \le k+2\\
0,& k+2 < |x|
\end{cases}\\
& \ge I_{[k,k+1]}(|x|)
\end{align*}
with the price $\mu_k = \nu_{k+2} - \nu_{k+1} - \nu_k + \nu_{k-1}$.
For $k=0$, $T_0(x)=1 - ((|x| - 1)_+ - (|x| - 2)_+)$, which is a single
trapezoid.  Symmetric bull spread and the positive side of symmetric
trapezoidal hedge are depicted in Figure \ref{fig:2} and Figure
\ref{fig:3}, respectively.

\begin{figure}[htbp]
\setlength{\unitlength}{5mm}
\begin{center}
\begin{picture}(20,3)(-10,0)
\put(-9,0){\line(1,0){18}}
\put(0,-0.2){\line(0,1){0.4}}
\put(-3,0){\line(-1,1){2}}
\put(-3.3,-1){$-k$}
\put(3,0){\line(1,1){2}}
\put(2.9,-1){$k$}
\put(-5,2){\line(-1,0){3}}
\put(-5,-0.1){\line(0,1){0.2}}
\put(5,-0.1){\line(0,1){0.2}}
\put(5,2){\line(1,0){3}}
\put(-6.5,-1){$-k-1$}
\put(4.3,-1){$k+1$}
\end{picture}
\end{center}
\caption{Symmetric bull spread}
\label{fig:2}
\end{figure}
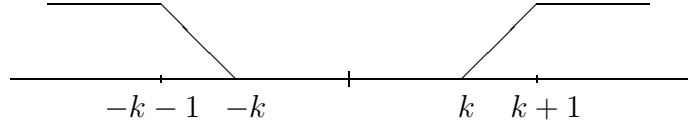

Now we state the following theorem.

\begin{theorem}
\label{thm:2}
Suppose that the set of hedges $\cH = \{ h_k(x)=(|x|-k)_+ \mid
k=0,1,2,\dots\}$ are available to Skeptic.
Then  in the unbounded forecasting game with $\cH$ 
Skeptic can force $\bar x_n \rightarrow 0$.
\end{theorem}

The rest of this section is devoted to a proof of this theorem.  As in
the previous section we prove it by a series of lemmas.

\begin{lemma} 
\label{lem:CBC1}
Under the condition of Theorem \ref{thm:2}
Skeptic can force $E_2$ in {\rm (\ref{eq:E2})}.
\end{lemma}

\begin{proof}  For $k\ge 1$, 
$
(|x| - k+1)_+ - (|x| - k)_+   \ge I_{[k,\infty)}(|x|)
$ and
\[
\sum_{n=1}^\infty \big((|x_n| - n+1)_+ - (|x_n| - n)_+ \big)  \ge
\sum_{n=1}^\infty I_{[n,\infty)}(|x_n|).
\]
The left-hand side can be bought with the total finite price of
$\nu_0$.
The rest of the proof is the same as in Lemma 
\ref{lem:BC0}.
\end{proof}

\begin{lemma}
Under the condition of Theorem \ref{thm:2}
Skeptic can force $E_3$ in {\rm (\ref{eq:E3})}.
\end{lemma}

\begin{proof}
  At round $n$ Skeptic is to buy $(k+1)^2$ units of the symmetric
  trapezoidal hedge $T_k$
  for each $k=0,1,\ldots,n-1$. We note
\[
\sum_{k=0}^{n-1} (k+1)^2 T_k(x_n) \ge 
x_n^2 I_{\{|x_n|\le n\} }.
\]
Dividing the above by $n^2$ and summing up over all rounds $n = 1, 2, \dots,$ we have
\[
\sum_{n=1}^\infty \sum_{k=0}^{n-1} \frac{(k+1)^2}{n^2} T_k(x_n) \ge
\sum_{n=1}^\infty \frac{x_n^2}{n^2} I_{\{|x_n|\le n\} }.
\]
Now we evaluate the total price of the left-hand side. Since $T_k$ is
available at each round, the price is the same if we replace $x_n$ by
$x_1$ in $T_k$. Then
\begin{align*}
\sum_{n=1}^\infty \sum_{k=0}^{n-1} \frac{(k+1)^2}{n^2}T_k(x_1) 
&= \sum_{n=1}^\infty \frac{1}{n^2} \sum_{l=1}^n l^2 T_{l-1} (x_1)
= \sum_{l=1}^\infty l^2 T_{l-1}(x_1) \sum_{n=l}^\infty \frac{1}{n^2} 
\\
&\le 2\sum_{l=1}^\infty l T_{l-1} (x_1)
\le 6|x_1|. 
\end{align*}
As noted above 
$|x_1|$ is available to Skeptic 
with finite price $\nu_0$, 
so that the left-hand side is also 
available to him with the total finite price 
\[
\sum_{n=1}^\infty \sum_{k=0}^{n-1} \frac{(k+1)^2}{n^2}\mu_k 
\le 6\nu_0.  
\]
The rest of the proof is the same as in Lemma \ref{lem:BC0}.
\end{proof}

In the following $x_n$  hedged by $(|x_n|-n)_+$ is denoted as
\[
x_{n,n}= x_n  +  (|x_n|-n)_+
=\begin{cases}
 -n, &  x_n < -n, \\
  x_n, & -n \le x_n \le n, \\
  2x_n - n, & x_n > n.
\end{cases}
\]
This has the price $\nu_n$.  Similarly we denote
$\tilde x_{n,n}=-x_n  +  (|x_n|-n)_+$ which is $-x_n$ hedged by 
$(|x_n|-n)_+$.  Note that
\[
x_{n,n} \ge -n, \qquad  \tilde x_{n,n} \ge -n.
\]
On $E_2$, $x_{n,n}$ and $x_n
I_{\{|x_n|\le n\}}$ differ only for finite number of $n$.  Therefore
conditional on $E_2$, Skeptic can force 
\begin{equation}
\label{eq:E3'}
E_3' = \{ \xi \mid \sum_n \frac{x_{n,n}^2}{n^2} < \infty\}, \qquad
E_3'{}' = \{ \xi \mid \sum_n \frac{\tilde x_{n,n}^2}{n^2} < \infty\}.
\end{equation}

\begin{lemma} Under the condition of Theorem \ref{thm:2} and
  conditional on $E_2$, Skeptic can force
\[
E_5 = \{ \xi \mid \sum_n \frac{(x_{n,n}-\nu_n)^2}{n^2} < \infty
\}
\]
\end{lemma}


\begin{proof} 
Since $(x_{n,n} - \nu_n)^2 \le 2 x_{n,n}^2 + 2 \nu_n^2$
\[
\sum_{n=1}^N \frac{(x_{n,n}-\nu_n)^2}{n^2} 
\le 
2\sum_{n=1}^N
\frac{x_{n,n}^2}{n^2} +  2\sum_{n=1}^N \frac{\nu_n^2}{n^2}
\le 2\sum_{n=1}^N
\frac{x_{n,n}^2}{n^2} +  2\nu_0^2 \frac{\pi^2}{6}
\]
By (\ref{eq:E3'}), conditional on $E_2$, Skeptic can force 
$\sum_{n=1}^\infty x_{n,n}^2/n^2 < \infty$.
Therefore conditional on $E_2$, he can force $E_5$
\end{proof}

Similarly Skeptic can force $E_5$ with $x_{n,n}$ replaced by $\tilde x_{n,n}$.

Finally the following lemma proves Theorem \ref{thm:2}
in conjunction with Kronecker's lemma.

\begin{lemma}  
Under the condition of Theorem \ref{thm:2}
Skeptic can weakly force 
\begin{align*}
E'_4 = \{ \xi \mid  \sum_n \frac{x_{n,n} - \nu_n}{n} \
\text{\rm converges to a finite value} \}
\end{align*}
conditional on $E_2 \cap E_3'$.
\end{lemma}

\begin{proof} 
We take $\epsilon$ as
\[
0 < \epsilon  < \frac{1}{2(1 + \nu_0)}, 
\]
and consider Skeptic's strategy betting $\epsilon
\cK_{n-1}/n$ on $x_{n,n}-\nu_n$ at round $n$. 
Then his capital at the end of round $n$ is
\[
{\cal K}_n = {\cal K}_{n-1}(1 + \frac{\epsilon}{n}(x_{n,n} - 
\nu_n)) 
= {\cal K}_0 \prod_{i=1}^n (1 + \frac{\epsilon}{i}(x_{i,i} - \nu_i)).   
\]
By the choice of $\epsilon$ and $|x_{i,i}/i|\le 1$,
\[
\frac{\epsilon}{i}(x_{i,i} - \nu_i) \ge - \frac{1}{2},
\]
so that from $\log (1+t) \ge t-t^2$ for $t \ge -1/2$, his log capital 
is bounded from below as
\begin{align*}
\log {\cal K}_n &\ge \log {\cal K}_0 + \epsilon \sum_{i=1}^n
\frac{x_{i,i}-\nu_i}{i}
 - \epsilon^2 \sum_{i=1}^n 
\frac{(x_{i,i} - \nu_i)^2}{i^2}.
\end{align*}
In the right-hand side the third term is bounded on $E_5$. 
By considering this inequality for the interval $n \le i \le n'$, we have 
\begin{align*}
\log {\cal K}_{n'} - \log {\cal K}_{n-1} 
\ge \epsilon \sum_{i=n}^{n'}
\frac{x_{i,i} - \nu_i}{i}
 - \epsilon^2 \sum_{i=n}^{n'} 
\frac{(x_{i,i} - \nu_i)^2}{i^2}.
\end{align*}
As in the proof of Lemma \ref{lem:e4}, considering both $x_{n,n}$ and 
$\tilde x_{n,n}$, 
$\log {\cal K}_n$ converges to a finite limit almost surely, 
and thus by Cauchy criterion we see that 
\[
\sum_{i=1}^n
\frac{x_{i,i} - \nu_i}{i}
\] 
converges almost surely. 
\end{proof}

As proved in Appendix \ref{app:2},
$\nu_n \rightarrow 0$ as $n\rightarrow\infty$.
Then by Kronecker's lemma we have
\[
\frac{1}{n}\sum_{i=1}^n (x_{i,i} - \nu_i) = 
\frac{1}{n}\sum_{i=1}^n x_{i,i} - \frac{1}{n}\sum_{i=1}^n \nu_i 
\to 0\quad \textrm{as}\ n \to \infty.
\]
In the above, 
$(1/n)\sum_{i=1}^n \nu_i  \rightarrow 0$ 
so that 
$\sum_{i=1}^n x_{i,i}/n$ also converges to 0. 
Since $x_n$ and $x_{n,n}$ differ only for finite number of $n$ on $E_2$, 
it is concluded that $\bar{x}_n$ converges to 0 almost surely. 
This completes the proof of Theorem \ref{thm:2}.

\section{Marcinkiewicz-Zygmund strong law}
\label{sec:MZ}

In this section we consider a remarkable generalization 
by Marcinkiewicz and Zygmund (See \cite{gut},  \cite{marcinkeiwicz-zygmund})
of Kolmogorov's measure-theoretic SLLN for i.i.d.\ random variables 
with finite expected value $E|x_n| < \infty$.
Marcinkiewicz-Zygmund strong law 
asserts that for i.i.d.\ random 
variables $x_1, x_2, \dots $ with $E|x_n|^r < \infty$ for $0 < r < 2$ 
and $Ex_n = 0$ when $1 \le r < 2$, the following measure-theoretic SLLN holds
\[
\frac{x_1 + \cdots + x_n}{n^{1/r}} \to 0\quad \textrm{as}\ n \to \infty \quad a.s.
\]

Considering the meaning of the hedge $|x|^r$ for $x$ in betting games, we
treat the case $1 < r < 2$ and for this case establish a
game-theoretic version of Marcinkiewicz-Zygmund SLLN. As noted in Proposition
\ref{prop:1}, the availability of a single $h(x) = |x|^r$ is again not
sufficient.  Hence here, for the game-theoretic Marcinkiewicz-Zygmund
SLLN we assume that the following set of hedges are available to
Skeptic.
\begin{equation}
\cH_r = \cH_{1r} \cup \cH_{2r}, \quad \text{where}\quad
\begin{cases}
\cH_{1r}=\{ h_{kr}(x)= (|x|^r-k)_+ \mid  k=0, 1, 2,\dots\} &\\
\cH_{2r}=\{ h_{k^{1/r}}(x)= (|x|-k^{1/r})_+ \mid  k=0, 1, 2,\dots\}&
\end{cases}.
\label{eq:countable-hedges-r}
\end{equation}
Let $\nu_{kr}$ denote the price of $h_{kr}(x)= (|x|^r-k)_+$ and
let $\nu_{k^{1/r}}$ denote the price of 
$h_{k^{1/r}}(x)= (|x|-k^{1/r})_+$.
Also 
assuming that Skeptic is allowed to sell and combine these hedges 
within his collateral duty, we state the following theorem.
\begin{theorem}
\label{thm:5}
Let $1 < r < 2$.  Suppose that the set of hedges $\cH_r$ in
{\rm (\ref{eq:countable-hedges-r})}
is available to Skeptic.
Then in the unbounded forecasting game with $\cH_r$ 
Skeptic can force $(x_1 + \cdots + x_n)/n^{1/r} \rightarrow 0$.
\end{theorem}

\begin{remark}
\label{ref:MZ-hedge}
In this theorem  $\cH_r$ consists of two sets of hedges $\cH_{1r}$ and
$\cH_{2r}$.  $\cH_{2r}$ is included in $\cH$ just for convenience.
Each $h_{k^{1/r}}(x)$ can be superreplicated and underreplicated by an
infinite combination of hedges from $\cH_{1r}$ and the theorem holds
without $\cH_{2r}$.  Since this makes the proof considerably
messier, we include $\cH_{2r}$ in the set of hedges.  We give more
discussion on this point in Section \ref{sec:discussions}.
\end{remark}

The proof of Theorem \ref{thm:5} proceeds almost in the same way as that of 
Theorem \ref{thm:2}.  However  we have to make different uses of
hedges from $\cH_{1r}$ and from $\cH_{2r}$.
At first we enumerate relevant events. 
\[
E_{2r} = \{ \xi \mid |x_n|^r \ge n   
\text{\ for only finite number of } n\}.
\]
\[
E_{3r} = \{ \xi \mid \sum_n \frac{x_n^2}{n^{2/r}}
I_{\{|x_n|^r\le n\}} < \infty \}.
\]

\begin{lemma} 
\label{lem:CBC1r}
Under the condition of Theorem \ref{thm:5}
Skeptic can force $E_{2r}$.
\end{lemma}
\begin{proof}  For $k\ge 1$, 
$
(|x|^r - (k-1))_+ - (|x|^r - k)_+   \ge I_{[k,\infty)}(|x|^r)
$ and
\[
\sum_{n=1}^\infty \big((|x_n|^r - (n-1))_+ - (|x_n|^r - n)_+ \big)  \ge
\sum_{n=1}^\infty I_{[n,\infty)}(|x_n|^r).
\]
The left-hand side can be bought with the total finite price of
$\nu_{0r}$.
The rest of the proof is the same as in Lemma 
\ref{lem:BC0}.
\end{proof}
\begin{lemma}
Under the condition of Theorem \ref{thm:5}
Skeptic can force $E_{3r}$.
\end{lemma}
\begin{proof}
Consider the following trapezoidal hedge
\begin{align*}
T_{kr}(x)& = (|x|^r - (k-1))_+ - (|x|^r - k)_+  - 
\big((|x|^r - (k+1))_+ - (|x|^r - (k+2))_+ \big)\\
&=
\begin{cases}
0,&   |x|^r \le k-1\\
|x|^r - (k-1), & k-1 < |x|^r \le  k\\
1, & k < |x|^r \le k+1\\
k+2 - |x|^r, & k+1 < |x|^r \le k+2\\
0,& k+2 < |x|^r
\end{cases}\\
&\ge I_{[k,k+1]}(|x|^r)
\end{align*}
with the price 
$\mu_{kr} = \nu_{k+2,r} - \nu_{k+1,r} -
\nu_{kr} + \nu_{k-1,r}$.

At round $n$ Skeptic is to buy $(k+1)^2$ units of the hedge $T_{kr}$ for each $k=0,1,\ldots,n-1$. We note 
\[
\sum_{k=0}^{n-1} (k+1)^2 T_{kr}(x_n) \ge 
x_n^2 I_{\{|x_n|^r\le n\} }.
\]
Dividing the above by $n^{2/r}$ and summing up over all rounds $n = 1, 2, \dots,$ we have
\[
\sum_{n=1}^\infty \sum_{k=0}^{n-1} \frac{(k+1)^2}{n^{2/r}} T_{kr}(x_n) \ge
\sum_{n=1}^\infty \frac{x_n^2}{n^{2/r}} I_{\{|x_n|^r\le n\} }.
\]
As in the previous section, for the consideration of the total price,
we can replace $T_{kr}(x_n)$ by $T_{kr}(x_1)$.
Then the left-hand side can be evaluated as 
{\allowdisplaybreaks
\begin{align*}
\sum_{n=1}^\infty \sum_{k=0}^{n-1} \frac{(k+1)^2}{n^{2/r}}T_{kr}(x_1) 
&= \sum_{n=1}^\infty \frac{1}{n^{2/r}} \sum_{l=1}^n l^2 T_{l-1,r}(x_1)  
= \sum_{l=1}^\infty l^2 T_{l-1,r}(x_1) \sum_{n=l}^\infty \frac{1}{n^{2/r}} \\
&\le \frac{2^{(2/r)-1}}{(2/r) - 1}
\sum_{l=1}^\infty \frac{1}{l^{(2/r)-1}}l^2 T_{l-1,r} (x_1)\\
&\le \frac{2^{(2/r)-1}}{(2/r) - 1}
\sum_{l=1}^\infty \frac{1}{l^{(2/r)-1}}(l^{1/r})^{2-r}l^r T_{l-1,r}(x_1)\\
&\le \frac{2^{(2/r)-1}}{(2/r) - 1}
\sum_{l=1}^\infty l^r T_{l-1,r}(x_1) \le \frac{3\cdot 2^{(2/r)-1}}{(2/r) - 1}|x_1|^r. 
\end{align*}
}
Since $|x_1|^r$ is available to Skeptic with finite price $\nu_{0r}$, 
the left-hand side is also available to him with the total finite price 
\[
\sum_{n=1}^\infty \sum_{k=0}^{n-1} \frac{(k+1)^2}{n^{2/r}}\mu_{kr} 
\le \frac{3\cdot 2^{(2/r)-1}}{(2/r) - 1}\nu_{0r}.
\]
\end{proof}

So far we have used hedges from $\cH_{1r}$ for forcing various events.
In the following $x_n$ will be hedged by elements from $\cH_{2r}$.
We hedge $x_n$ by $h_{n^{1/r}}(x_n)=(|x_n|-n^{1/r})_+$.
Write
\[
x_{nn,r}
= x_n + (|x_n| - n^{1/r})_+ .
\]
This has the price $\nu_{n^{1/r}}$. On $E_{2r}$, $x_{nn,r}$ and $x_n
I_{\{|x_n|^r\le n\}}$ differ only for finite number of $n$.  Therefore
conditional on $E_{2r}$, Skeptic can force 
\begin{equation}
\label{eq:E3'r}
E_{3r}' = \{ \xi \mid \sum_n \frac{x_{nn,r}^2}{n^{2/r}} < \infty\}.
\end{equation}

\begin{lemma} Under the condition of Theorem \ref{thm:5} and
  conditional on $E_{2r}$, Skeptic can force
\[
E_{5r} = \{ \xi \mid \sum_n \frac{(x_{nn,r}-\nu_{n^{1/r}})^2}
{n^{2/r}} < \infty\}.
\]
\end{lemma}

\begin{proof} 
\[
\sum_n \frac{(x_{nn,r}-\nu_{n^{1/r}})^2}{n^{2/r}} \le 2\sum_n  
  \frac{x_{nn,r}^2}{n^{2/r}}
+ 2\sum_n \frac{\nu_{n^{1/r}}^2}{n^{2/r}}.
\]
Both terms are finite on $E_{2r}$.
\end{proof}

Now we use the $\epsilon$-strategy as before.

\begin{lemma}  
Under the condition of Theorem \ref{thm:5}
Skeptic can weakly force 
\begin{align*}
E'_{4r} = \{ \xi \mid  \sum_n \frac{x_{nn,r} - \nu_{n^{1/r}}}{n^{1/r}} \ 
\text{\rm converges to a finite value} \}
\end{align*}
conditional on $E_{2r} \cap E_{3r}'$.
\end{lemma}

\begin{proof} 
We take $\epsilon$ as
\[
0 < \epsilon < \frac{1}{2(1 + \nu_{0})}, 
\]
and consider Skeptic's strategy betting $\epsilon
\cK_{n-1}/n$ on $x_{n,n}-\nu_{n^{1/r}}$ at round $n$. 
Then his capital at the end of round $n$ is
\[
{\cal K}_n = {\cal K}_{n-1}(1 + \frac{\epsilon}{n}(x_{nn,r} - 
\nu_{n^{1/r}})) 
= {\cal K}_0 \prod_{i=1}^n (1 + \frac{\epsilon}{i}(x_{ii,r} - 
\nu_{i^{1/r}})).   
\]
By the choice of $\epsilon$ and $|x_{ii,r}/i^{1/r}|\le 1$,
\[
\frac{\epsilon}{i^{1/r}}(x_{ii,r} - \nu_{i^{1/r}}) \ge - \frac{1}{2},
\]
so that from $\log (1+t) \ge t-t^2$ for $t \ge -1/2$, his log capital 
is bounded from below as
\begin{align*}
\log {\cal K}_n &\ge \log {\cal K}_0 + \epsilon \sum_{i=1}^n
\frac{x_{ii,r} - \nu_{i^{1/r}}}{i^{1/r}}
 - \epsilon^2 \sum_{i=1}^n 
\frac{(x_{ii,r} - \nu_{i^{1/r}})^2}{i^{2/r}}.
\end{align*}
In the right-hand side the third term is bounded on $E_{5r}$.  
By considering this inequality for the interval $n \le i \le n'$,  
\begin{align*}
\log {\cal K}_{n'} - \log {\cal K}_{n-1} 
\ge \epsilon \sum_{i=n}^{n'}
\frac{x_{ii,r} - \nu_{i^{1/r}}}{i^{1/r}}
 - \epsilon^2 \sum_{i=n}^{n'} 
\frac{(x_{ii,r} - \nu_{i^{1/r}})^2}{i^{2/r}}.
\end{align*}
In the above $\log {\cal K}_n$ converges to a finite limit almost surely, 
and thus as before  
\[
\sum_{i=1}^n
\frac{x_{ii,r} - \nu_{i^{1/r}}}{i^{1/r}}
\] 
converges almost surely. 
\end{proof}

We now need to take care of  $n^{-1/r} \sum_{i=1}^n \nu_{i^{1/r}}$.
\begin{lemma}
\label{lem:cesaro-r}
\[
\frac{1}{n^{1/r}} \sum_{i=1}^n \nu_{i^{1/r}} \rightarrow 0
\qquad \text{\rm as} \quad n\rightarrow\infty.
\]
\end{lemma}

\begin{proof}
Since $r>1$ 
\[
\big( 1- \frac{n^{1/r}}{|x|}\big)_+ \le \big (1-
\frac{n}{|x|^r}\big)_+ , \qquad \forall x.
\]
Also for $|x| \ge n^{1/r}$ we have $|x|^{r-1} \ge n^{1-1/r}$. Therefore
\begin{align*}
(|x| - n^{1/r})_+ 
&=  |x| \big( 1- \frac{n^{1/r}}{|x|}\big)_+ \\
&\le \frac{|x|^{r-1}}{n^{1-1/r}} |x|\big( 1- \frac{n}{|x|^r}\big)_+ \\
& = n^{1/r-1} (|x|^r - n)_+
\end{align*}
It follows that the prices of $(|x| - n^{1/r})_+ $ and 
$(|x|^r - n)_+$ have to satisfy
\[
\nu_{n^{1/r}} \le n^{1/r-1} \nu_{nr}.
\]
Therefore 
\[
\frac{1}{n^{1/r}} \sum_{i=1}^n \nu_{i^{1/r}} \le 
\frac{1}{n^{1/r}} \sum_{i=1}^n i^{1/r-1} \nu_{ir}.
\]
Note that $\nu_{ir}\rightarrow 0$ as $i \rightarrow \infty$ by the argument in
Appendix \ref{app:2}.
Then the right-hand side converges to 0 as $n\rightarrow\infty$ by 
Ces\`aro's lemma (12.6 of \cite{williams}).
\end{proof}

Now by an extended form of Kronecker's lemma (12.7 of
\cite{williams})
\[
\frac{1}{n^{1/r}}\sum_{i=1}^n (x_{ii,r} - \nu_{i^{1/r}}) = 
\frac{1}{n^{1/r}}\sum_{i=1}^n x_{ii,r} - 
\frac{1}{n^{1/r}}\sum_{i=1}^n \nu_{i^{1/r}} 
\to 0\quad \textrm{as}\ n \to \infty.
\]
so that 
$\sum_{i=1}^n x_{i,i}/n^{1/r}$ also converges to 0. 
Since $x_n$ and $x_{nn,r}$ differ only for finite number of $n$ on  $E_{2r}$, 
it follows that $(x_1 + \cdots + x_n)/n^{1/r}$ converges to 0 almost surely. 
This completes the proof of Theorem \ref{thm:5}.

\section{Some discussions}
\label{sec:discussions}

In this paper we proved various game-theoretic versions of SLLN for
unbounded variables.  In Section \ref{sec:with-truncation} we proved a
version corresponding to the sample average of i.i.d.\
measure-theoretic random variables.  There we assumed availability of
countable symmetric call option type hedges.  We chose this set of
hedges for convenience and concreteness.  Other choices are equally
conceivable, as long as the set of hedges is rich enough to produce
step-function type hedges (cf. Figure \ref{fig:2}).  

We might as well assume that if a hedge $h$ is available to Skeptic,
all other hedges weaker then $h$ are available to him with price no
more than that of $h$.  We call a set of hedges $\cH$ {\rm weakly
  closed} if
\[
h \in \cH,  \ 0 \le g(x) \le h(x), \forall x\in {\mathbb R} 
\quad  \Rightarrow \quad  g\in \cH.
\]
We might argue that this is a reasonable assumption, because if $h$ is
available to Skeptic, he can ask to buy a weaker $g$ with the same
price as $h$ and someone should be willing to sell $g$ to Skeptic with
the same price, because it presents an arbitrage opportunity to the
seller.  If $\cH$ is weakly closed, then for each $t\in {\mathbb R}$
\[
I_{(-\infty,t]}(x)
\]
has to be available to Skeptic.  This shows that if $\cH$ is weakly
closed, then the entire distribution function of the Reality's move
$x$ is priced in the game.  The assumption of weakly closed $\cH$
seems to be too strong from game-theoretic viewpoints.  However we
should mention that in measure-theoretic proofs the probability distribution is
assumed and truncation is freely used.

The discussion on generality of probability games in Chapter 8 of
Shafer and Vovk (2001) convincingly argues that measure-theoretic
martingales can be reduced to game-theoretic martingales.  If we
interpret Theorem \ref{thm:1} in measure-theoretic terms and just
rewrite our proof in measure-theoretic terms, we obtain the following
result.

\begin{proposition}
\label{prop:app}
Let $\{Y_n\}$ be a measure-theoretic martingale adapted to an
increasing family of $\sigma$-fields $\{{\cal F}_n\}$.  Let $h$ be a
function satisfying {\rm (A1)--(A3)}.  If the measure-theoretic
conditional expectation
\[
E[h(Y_n - Y_{n-1}) \mid {\cal F}_{n-1}]
\]
is uniformly bounded, then $P(\lim_n Y_n/n=0)=1$.
\end{proposition}
Except for Proposition \ref{prop:1} we could avoid measure theory to
establish our theorems.  We believe that this again shows 
effectiveness of game-theoretic proofs  as we have shown in our previous works
(\cite{ktt1}, \cite{kumon-takemura}).

For the Marcinkiewicz-Zygmund strong law in Section \ref{sec:MZ} we
have given a game-theoretic proof for $r>1$.  We also assumed
availability of two kinds of hedges for convenience as we discussed in
Remark \ref{ref:MZ-hedge}.  If we make the blanket assumption that
$\cH$ is weakly closed, then we believe that measure-theoretic proof
of the Marcinkiewicz-Zygmund strong for $0 < r < 1$ can be translated
to game-theoretic proof without too many modifications.  From
game-theoretic viewpoint however, the case $r < 1$ does not seem to be
natural.


\appendix
\section{Proofs of Proposition \ref{prop:1} and Proposition \ref{prop:a3}}
\label{sec:app1}

\noindent
{\it Proof of Proposition \ref{prop:1}.} \qquad
We argue by contradiction. Suppose there exists Skeptic's strategy
$\cP$ which allows Skeptic to observe his collateral duty with the
initial capital $\cK_0=1$ and $\lim_n\cK^\cP_n  = \infty$ whenever
$s_n/n^{1/r} \not\rightarrow 0$, where $s_n = x_1 + \cdots + x_n$.
Consider a random strategy of Reality, where each $x_n$, $n > \nu$, is
independently chosen as
\[
P(x_n=0)= 1 - \frac{\nu}{n}, \ P(x_n =n^{1/r})=
P(x_n=- n^{1/r})=\frac{\nu}{2n}.
\]
Here $\nu$ is the price of $h(x)=|x|^r$. 
Then by the second part of measure-theoretic Borel-Cantelli lemma
\begin{equation}
\label{eq:a1}
1 = P(|x_n|=n^{1/r} \ \ i.o.)= P(|x_n|/n^{1/r}=1 \ \ i.o.).
\end{equation}
Note that if $s_n/n^{1/r} \rightarrow 0$, then 
$x_n/n^{1/r} \rightarrow 0$ because
\[
\frac{s_n}{n^{1/r}} = \Big(\frac{n-1}{n}\Big)^{1/r}
\frac{s_{n-1}}{(n - 1)^{1/r}} + \frac{x_n}{n^{1/r}}.
\]
Therefore (\ref{eq:a1}) implies that $P(s_n/n^{1/r} \rightarrow 0)=0$.
Then by our assumption $P(\cK^\cP_n \rightarrow\infty)=1$.  However
under the randomized strategy of Reality 
$\cK^\cP_n$ is a measure-theoretic non-negative martingale and its 
measure-theoretic expectation is  $E(\cK^\cP_n)=\cK_0=1$.  Then by
Doob's martingale inequality (e.g.\ Theorem 14.6 of \cite{williams})
\[
P(\max_{k \le n} \cK^\cP_k \ge c)\le \frac{1}{c}, \quad \forall c>0, \forall n.
\]
and $P(\sup_n \cK^\cP_n \ge c)\le 1/c$. 
But this contradicts $P(\cK^\cP_n \rightarrow\infty)=1$. \hfill $\Box$.

\bigskip
\noindent
{\it Proof of Proposition \ref{prop:a3}.} \qquad 
Consider a random strategy of Reality, where each $x_n$ for $n$,
$h(n) > \nu$,  is independently chosen as 
\[
P(x_n=0)= 1 - \frac{\nu}{h(n)}, \ P(x_n =n)=
P(x_n=-n)=\frac{\nu}{2h(n)}.
\]
The rest of the proof is the same as the proof of Proposition 
\ref{prop:1}. \hfill $\Box$

\section{Proof of the fact $\lim_{k\rightarrow\infty} \nu_k=0$}
\label{app:2}

Consider the identity for $x\in {\mathbb R}$:
\begin{equation}
\label{eq:arbitrage0}
|x| = \sum_{k=0}^\infty \big( (|x|-k)_+ - (|x|-k-1)_+ \Big).
\end{equation}
For each real $x$, the right-hand side is actually a  finite sum and
there is no question on the convergence.  On the other hand consider
the identity
\[
\nu_0 = \sum_{k=0}^{K-1} (\nu_k - \nu_{k+1}) + \nu_K.
\]
Since $\{\nu_k\}$ is a monotone non-increasing sequence of non-negative reals
\[
c = \lim_{K\rightarrow\infty} \nu_K  \ge 0
\]
exists. If $c>0$ then, $\nu_0 > \sum_{k=0}^{\infty} (\nu_k -
\nu_{k+1})$.  But then Skeptic can sell $|x|$ and buy the right-hand
side of (\ref{eq:arbitrage0}) and he is certain to make money. This
contradicts the assumption of coherence.

\end{document}